\documentclass{icm2010}

\title[Locally homogeneous geometric manifolds]
{Locally homogeneous geometric manifolds
}
\author[W. Goldman]{William M. Goldman}

\contact[wmg@math.umd.edu]{Department of Mathematics \\
University of Maryland \\ College Park, MD 20742 USA}

\newtheorem{theorem}{Theorem}
 \newtheorem*{theorem*}{Theorem}

\begin{document}

\thanks{Partially supported by the National Science Foundation}

\newcommand{\R}{\mathbb{R}}
\newcommand{\E}{\mathbb{E}}
\newcommand{\Z}{\mathbb{Z}}
\newcommand{\C}{\mathbb{C}}
\newcommand{\rpt}{{\R\mathsf{P}^2}}
\newcommand{\rpn}{{\R\mathsf{P}^n}}
\newcommand{\rpthree}{{\R\mathsf{P}^3}}
\newcommand{\Hom}{{\mathsf{Hom}}}
\newcommand{\Ham}{{\mathsf{Ham}}}
\newcommand{\tM}{\tilde{M}}
\newcommand{\dev}{\mathsf{dev}}
\newcommand{\tr}{\mathsf{tr}}
\newcommand{\Dev}{\mathsf{Dev}}
\newcommand{\Ff}{\mathcal{F}}
\newcommand{\PSLtZ}{\mathsf{PSL}(2,\Z)}
\newcommand{\PSLtR}{\mathsf{PSL}(2,\R)}
\newcommand{\SLtZ}{\mathsf{SL}(2,\Z)}
\newcommand{\Inn}{\mathsf{Inn}}
\newcommand{\Out}{\mathsf{Out}}
\newcommand{\Aut}{\mathsf{Aut}}
\newcommand{\hol}{\mathsf{hol}}
\newcommand{\Euler}{\mathsf{Euler}}
\newcommand{\cpo}{\C\mathsf{P}^1}
\newcommand{\PSLtC}{\mathsf{PSL}(2,\C)}
\newcommand{\SLtC}{\mathsf{SL}(2,\C)}
\newcommand{\SLtR}{\mathsf{SL}(2,\R)}
\newcommand{\PGLnoR}{\mathsf{PGL}(n+1,\R)}
\newcommand{\Spfour}{\mathsf{Sp}(4,\R)}
\newcommand{\SLfourR}{\mathsf{SL}(4,\R)}
\newcommand{\PUno}{\mathsf{PU}(n,1)}
\newcommand{\Spnr}{\mathsf{Sp}(2n,\R)}
\newcommand{\Ad}{\mathsf{Ad}}
\newcommand{\SLthreeR}{\mathsf{SL}(3,\R)}

\newcommand{\DGXS}{\mathsf{Def}^{(G,X)}(\Sigma)}
\newcommand{\DpGXS}{\mathsf{Def}^{(G,X)}(\Sigma)}
\newcommand{\Mod}{\mathsf{Mod}(\Sigma)}

\newcommand{\Fricke}{\mathfrak{F}(\Sigma)}
\newcommand{\Teich}{\mathfrak{T}(\Sigma)}
\newcommand{\Cc}{\mathfrak{C}(\Sigma)}
\newcommand{\graph}{\mathsf{graph}}
\newcommand{\Oto}{\mathsf{O}(2,1)}
\newcommand{\Bb}{\mathbb{B}}

\newcommand{\HompGG}{\Hom\big(\pi_1(\Sigma),G\big)/G}
\newcommand{\HompG}{\Hom(\pi_1(\Sigma),G)/G}
\newcommand{\Zo}{\mathsf{Z}^1(\Sigma,\mathfrak{g}_{\Ad \rho})}

\newcommand{\MR}{Math.\ Rev.\ }

\begin{abstract}

Motivated by Felix Klein's notion that geometry is governed by its group of symmetry
transformations, Charles Ehresmann initiated the study of geometric structures on topological
spaces locally modeled on a homogeneous space of a Lie group. These locally homogeneous spaces
later formed the context of Thurston's $3$-dimensional geometrization program. The basic problem is for
a given topology $\Sigma$ and a geometry $X = G/H$, to classify all the possible ways of introducing the local geometry
of $X$ into $\Sigma$. 
For example, a sphere admits no local Euclidean geometry: there is no metrically accurate
Euclidean atlas of the earth. One develops a space whose points are equivalence classes of geometric structures
on $\Sigma$, which itself exhibits a rich geometry and symmetries arising from the topological symmetries of $\Sigma$.

We survey several examples of the classification of locally homogeneous geometric structures on manifolds in low dimension,
and how it leads to a general study of surface group representations.
In particular geometric structures are a useful tool in understanding
local and global properties of deformation spaces of representations of fundamental groups.
\end{abstract}



\maketitle


\section{Historical background}

While geometry involves quantitative measurements
and rigid metric relations, topology deals  with the
loose quantitative organization of points. Felix Klein proposed
in his 1872 Erlangen Program that the classical
geometries be considered as the properties of a 
space invariant under a transitive Lie group action.
Therefore one may ask which topologies support a system
of local coordinates modeled on a fixed homogeneous space
$ X = G/H$ such that on overlapping coordinate patches, the
coordinate changes  are locally restrictions of transformations
from $G$. 

In this generality this question was first asked by
Charles Ehresmann~\cite{Ehresmann1} 
at the conference ``Quelques questions de Geom\'etrie et de Topologie,'' in Geneva in 1935.
Forty years later, the subject of such {\em locally homogeneous geometric structures\/} experienced
a resurgence when W.\ Thurston placed his 3-dimensional geometrization program~\cite{Thurston_1979notes}
in the context of locally homogeneous (Riemannian) structures. The rich 
diversity of geometries on homogeneous spaces brings in a wide range of techniques, and the field has thrived through their interaction.

Before Ehresmann, the subject may be traced to several independent threads in the 19th century:
\begin{itemize}
\item The theory of monodromy of Schwarzian differential equations on Riemann surfaces, which arose from the integration of algebraic functions;
\item Symmetries of crystals led to the enumeration (1891) by Fedorov,
Sch\"oenflies and Barlow of the $230$  three-dimensional crystallographic
{\em space groups\/} (the $17$ two-dimensional {\em wallpaper groups\/} had been known much earlier). The general qualitative
classification of crystallographic groups is due to Bieberbach.
\item The theory of connections, curvature and parallel transport in Riemannian geometry, which arose from the classical theory of surfaces
in $\R^3$.
\end{itemize}
The uniformization of Riemann surfaces linked complex analysis to Euclidean and non-Euclidean
geometry. Klein, Poincar\'e and others saw that the moduli of Riemann surfaces, first conceived
by Riemann, related (via uniformization) to the deformation theory of geometric structures.
This in turn related to deforming discrete groups (or more accurately, representations of fundamental
groups in Lie groups), the viewpoint of the text of Fricke-Klein~\cite{FrickeKlein}.

\section
{The Classification Question\/}
Here is the fundamental general problem: 
Suppose we are given a manifold
$\Sigma$ (a {\em topology\/}) 
and a homogeneous space $(G, X = G/H)$ (a {\em geometry\/}).
Identify a space whose points correspond to equivalence
classes of $(G,X)$-structures on $\Sigma$. 
This space should inherit an action of the group of topological symmetries ({\em the mapping class group\/} $\Mod$) of 
$\Sigma$.
That is, how many inequivalent ways can one weave the geometry
of $X$ into the topology of $\Sigma$?
Identify the natural $\Mod$-invariant geometries on
this deformation space.

\section
{Ehresmann structures and development\/}
For $n>1$, the sphere $S^n$ admits no Euclidean structure.
This is just the familiar fact there is no metrically accurate atlas
of the world.
Thus the deformation space of Euclidean structures on $S^n$
is empty. On the other hand, the torus admits a rich class of
Euclidean structures, and (after some simple normalizations)
the space of Euclidean structures on $T^2$ identifies with the
quotient of the upper half-plane $H^2$ by the modular group
$\mathsf{PGL}(2,\Z)$. 

Globalizing the coordinate charts in terms of the 
{\em developing map\/} is useful here.
Replace the coordinate atlas by
a universal covering space $\tM\longrightarrow M$ 
with covering group $\pi_1(M)$.
Replace the coordinate charts by a local diffeomorphism,
the {\em developing map\/}  $\tM\xrightarrow{\dev} X$, 
as follows.
$\dev$ is equivariant with respect
to the actions of $\pi_1(M)$ by deck transformations on $\tM$
and by a representation $\pi_1(M)\xrightarrow{h} G$, respectively. 
The coordinate changes are replaced by the {\em holonomy
homomorphism\/} $h$. 
The resulting {\em developing pair\/} $(\dev,h)$ is unique up
to composition/conjugation by elements in $G$. This determines
the structure. 

Here is the precise correspondence.
Suppose that 
\begin{equation*}
\{(U_\alpha,\psi_\alpha)\mid U_\alpha\in \mathcal{U}\}
\end{equation*}
is a {\em $(G,X)$-coordinate atlas:\/}
$\mathcal{U}$ is an open covering by coordinate patches $U_\alpha$,
with coordinate charts $U_\alpha \xrightarrow{\psi_\alpha} X$
for $U_\alpha\in\mathcal{U}$. For every nonempty connected open subset $U\subset U_\alpha\cap U_\beta$, there is a (necessarily
unique) 
\begin{equation*}
g(U;U_\alpha,U_\beta)\in G
\end{equation*}
such that 
\begin{equation*}
\psi_\alpha|_{U} \;=\; g(U) \circ \psi_\beta|_{U}.
\end{equation*}

(Since a homogeneous space $X$ carries a natural real-analytic
structure invariant under $G$, every $(G,X)$-manifold carries
an underlying real-analytic structure. 
For convenience, therefore, we fix a smooth structure on 
$\Sigma$, and work in the differentiable category,
where tools such as transversality are available.
Since we concentrate here in low dimensions (like  $2$),
restricting to smooth manifolds and mappings sacrifices no
generality. Therefore, when we speak of ``a topological space $\Sigma$"
we really mean a smooth manifold $\Sigma$ rather than just
a topological space.)

The coordinate changes $\{g(U;U_\alpha,U_\beta)\}$ define
a {\em flat $(G,X)$-bundle\/} as follows.
Start with the trivial $(G,X)$-bundle
over the disjoint union $\coprod _{U_\alpha\in\mathcal{U}}  U_\alpha$,
having components
\begin{equation*}
E_\alpha :=  U_\alpha \times X  \xrightarrow{\Pi_\alpha} U_\alpha.
\end{equation*}
Now identify, for 
\begin{equation*}
(u,u_\alpha,u_\beta)\in U\times U_\alpha \times U_\beta,
\end{equation*} 
the two local total spaces
$U\times X\; \subset\; E_\alpha$ with 
$U\times X \;\subset\; E_\beta$ by
\begin{equation}\label{eq:identification}
\big(u,x\big)_\alpha 
\;\longleftrightarrow\; 
\big(u, g(U;U_\alpha,U_\beta) x\big)_\beta.
\end{equation}
The fibrations $\Pi_\alpha$ over $U_\alpha$ piece together to form
a fibration 
$E(M) \xrightarrow{\Pi} M$
over $M$ with fiber $X$, and structure group $G$,
whose total space $E = E(M)$ is the quotient space
of the $E_\alpha$ by the identifications \eqref{eq:identification}.
The foliations $\Ff_\alpha$ of $E_\alpha$ defined locally by the projections $U_\alpha\times X\longrightarrow X$ piece together to define a foliation  $\Ff(M)$ of $E(M)$ transverse to the fibration.
In this atlas, the coordinate changes are locally constant maps
$U_\alpha\cap U_\beta \longrightarrow G$. 
This {\em reduces the structure group\/} from $G$ with its manifold
topology to $G$ {\em with the discrete topology.}
We call the fiber bundle $\big(E(M),\Ff(M)\big)$ the {\em flat $(G,X)$-bundle tangent to $M$.}

Such a bundle pulls back to a trivial bundle over the universal
covering $\tM\longrightarrow M$. Thus it may be reconstructed
from the trivial bundle $\tM\times X \longrightarrow \tM$ as the
quotient of a $\pi_1(M)$-action on $\tM\times X$ covering the
action on $\tM$ by deck transformations. Such an action is
determined by a homomorphism $\pi_1(M)\xrightarrow{h} G$,
the {\em holonomy representation.}  Isomorphism classes of flat
bundles with structure group $G$ correspond to $G$-orbits
on $\Hom\big(\pi_1(M),G\big)$ by left-composition with inner automorphisms of $G$.

The coordinate charts $U_\alpha\xrightarrow{\psi_\alpha} X$
globalize to a section of the flat $(G,X)$-bundle $E\longrightarrow M$
as follows. The graph $\graph(\psi_\alpha)$ is a section
transverse both to the fibration and the foliation $\Ff_\alpha$.
Furthermore the identifications \eqref{eq:identification} imply that
the restrictions of $\graph(\psi_\alpha)$  and
$\graph(\psi_\beta)$ to $U\subset U_\alpha\cap U_\beta$
identify. Therefore all the $\psi_\alpha$ are the restrictions of a globally
defined $\Ff$-transverse section $M \xrightarrow{\Dev} E$.
We call this section the {\em developing section\/} since it exactly
corresponds to a developing map.

Conversely, suppose that $(E,\Ff)$ is a flat $(G,X)$-bundle over $M$ and $M\xrightarrow{s}E$ is a section transverse to $\Ff$.
For each $m\in M$, choose an open neighborhood $U$ such that
the foliation $\Ff$ on the local total space $\Pi^{-1}(U)$ is defined
by a submersion $\Pi^{-1}(U)\xrightarrow{\Psi_U} X$. 
Then the compositions $\Psi_U\circ s$ define  coordinate charts
for a $(G,X)$-structure on $M$.

In terms of the universal covering space $\tM\longrightarrow M$ 
and holonomy representation $h$, a section $M\xrightarrow{s}E$ corresponds to an $\pi_1(M)$-equivariant mapping 
$\tM \xrightarrow{\tilde{s}} X$, where $\pi_1(M)$ acts on $X$ via $h$.
The section $s$ is transverse to $\Ff$ if and only if the corresponding
equivariant map $\tilde{s}$ is a local diffeomorphism.

\section{Elementary consequences}
As the universal covering $\tilde{M}$ immerses in $X$, 
no $(G,X)$-structure exists when $M$ is closed with finite
fundamental group and $X$ is noncompact.
Furthermore if $X$ is compact and simply connected, then every
closed $(G,X)$-manifold with finite fundamental group would be
a quotient of $X$. 
Thus by extremely elementary considerations, 
no counterexample to the Poincar\'e conjecture could be
modeled on $S^3$.

When $G$ acts properly on $X$ (that is, when the isotropy group
is compact), then $G$ preserves a Riemannian metric on $X$
which passes down to a metric on $M$. This metric lifts to
a Riemannian metric on the the universal covering $\tilde{M}$,
for which $\dev$ is a local isometry.
Suppose that $M$ is closed. The Riemannian metric on $M$
makes $M$ into a metric space, which is necessarily complete.
By the Hopf-Rinow theorem,  $M$ is geodesically complete,
and (after possibly  replacing $X$ with its universal covering space
$\tilde{X}$, and $G$ by an appropriate group $\tilde{G}$ of lifts),
the local isometry $\dev$ is a covering space, and maps $\tM$ bijectively to $\tilde{X}$.
In particular such structures correspond to discrete cocompact
subgroups of $\tilde{G}$. 
In this way the subject of Ehresmann geometric structures
extends the subject of discrete subgroups of Lie groups.

In general, even for closed manifolds, 
the developing map may fail to be surjective (for example,
Hopf manifolds), and even may not be a covering space onto
its image (Hejhal~\cite{Hejhal}, Smillie~\cite{Smillie_thesis},
Sullivan-Thurston~\cite{SullivanThurston}).

\section{The hierarchy of geometries}

Often one geometry ``contains'' another geometry as follows.
Suppose that $G$ and $G'$ act transitively on $X$ and $X'$
respectively, and $X\xrightarrow{f}X'$ is a local diffeomorphism
equivariant respecting a homomorphism $G\xrightarrow{F}G'$.
Then (by composition with $f$ and $F$) every $(G,X)$-structure
determines a $(G',X')$-structure. For example, when $f$ is the
identity, then $G$ may be the subgroup of $G'$ preserving
some extra structure on $X=X'$. In this way, various flat pseudo-Riemannian geometries are refinements of affine geometry.
The three constant curvature Riemannian geometries (Euclidean,
spherical, and hyperbolic) have both realizations in conformal
geometry of $S^n$ (the Poincar\'e model) and in projective geometry
(the Beltrami-Klein model) in $\rpn$. In more classical differential-geometric terms, this is just the fact that the constant curvature Riemannian geometries are {\em conformally flat\/} (respectively
{\em projectively flat\/)}. Identifying conformal classes of conformally flat
Riemannian metrics as  Ehresmann structures follows from
Liouville's theorem on the classification of conformal maps of
domains in $\R^n$ for $n\ge 3$.

An interesting and nontrivial example is the classification of
closed similarity manifolds by Fried~\cite{Fried_sim}.
Here $X=\R^n$ and $G$ is its group of similarity transformations.
Fried showed that every closed $(G,X)$-manifold $M$ is either a 
Euclidean manifold (so $G$ reduces to the group of {\em isometries\/})
or a {\em Hopf manifold,\/}  a quotient of $\R^n\setminus\{0\}$
by a cyclic group of linear expansions. In the latter case $M$
carries a $\big(\R^+\cdot\mathsf{O}(n),\R^n\setminus\{0\}\big)$-structure.
Such manifolds are finite quotients of $S^{n-1}\times S^1$.
\section
{Deforming Ehresmann structures\/}

One would like a space whose points are equivalence classes
of $(G,X)$-structures on a fixed topology $\Sigma$.
The prototype of such a {\em deformation space\/} is the
{\em Teichm\"uller space\/} $\Teich$  of biholomorphism
classes of complex structures on a fixed surface $\Sigma$.
That is, we consider a Riemann surface $M$ with a 
diffeomorphism $\Sigma\longrightarrow M$, which is commonly
called a {\em marking.\/}
Although complex structures are not Ehresmann structures,
there is still a formal similarity. (This formal similarity can be
made into an equivalence of categories via the uniformization 
theorem, but this is considerably deeper than the present
discussion.) For example, every Riemann surface diffeomorphic to $T^2$ arises as $\C/\Lambda$, where $\Lambda\subset \C$ is
a lattice. Two such lattices $\Lambda, \Lambda'$ determine
isomorphic Riemann surfaces if $\exists\zeta\in\C^*$ such that 
$\Lambda' = \zeta\Lambda$. The space of such equivalence classes
identifies with the quotient $H^2/\PSLtZ$. 
The quotient $H^2/\PSLtZ$ 
has the natural structure of an {\em orbifold,\/}) and is not naturally a manifold.

In general deformation spaces will have
very bad separation properties. (For example the space of complete
affine structures on $T^2$ naturally identifies with the quotient of
$\R^2$ by the usual linear action of $\SLtZ$ (Baues, see
\cite{BauesGoldman}.) This quotient admits
no nonconstant continuous mappings into any Hausdorff space!)
To deal with such pathologies, we form a {\em larger space\/} with a group action, whose orbit space parametrizes isomorphism classes
of $(G,X)$-manifolds diffeomorphic to $\Sigma.$ 
In general, passing to the orbit space alone loses too much information,
and may result in an unwieldy topological space.
For this reason, considering the 
{\em deformation groupoid, \/}  consisting of structures
(rather than equivalence classes) and isomorphisms between
them, is a more meaningful and useful object to parametrize geometric structures.

Therefore we fix a smooth manifold $\Sigma$ and 
define a {\em marked $(G,X)$-structure\/} on $\Sigma$
as a pair $(M,f)$ where $M$ is a $(G,X)$-manifold and
$\Sigma\xrightarrow{f}M$  a diffeomorphism.
Suppose that $\Sigma$ is compact (possibly $\partial\Sigma\neq\emptyset$). 
Fix a fiber bundle $E$ over $\Sigma$ with fiber $X$ and structure
group $G$.
Give the set $\DpGXS$ of such marked $(G,X)$-structures on $\Sigma$
the $C^1$-topology on pairs $(\Ff,\Dev)$ of foliations $\Ff$
and smooth sections $\Dev$.
Clearly the diffeomorphism  group $\mathsf{Diff}(\Sigma)$ acts
on $\DpGXS$ by left-composition.
Define marked $(G,X)$-structures 
$(M,f)$ and $(M',f')$  to be {\em isotopic\/} if they are related
by an diffeomorphism of $\Sigma$ isotopic to the identity.

Define the {\em deformation space\/} of isotopy classes of marked 
$(G,X)$-structures on $\Sigma$ as the quotient space
\begin{equation*}
\DGXS \;:=\;  \DpGXS / \mathsf{Diff}_0(\Sigma).
\end{equation*}
Clearly the {\em diffeotopy group\/} 
$\pi_0\big(\mathsf{Diff}(\Sigma)\big)$
(which for compact surfaces $\Sigma$ is the 
{\em mapping class group\/} $\mathsf{Mod}(\Sigma)$
acts on the deformation space.

\section{Representations of the fundamental group}

The set of isomorphism classes of 
flat $G$-bundles over $\Sigma$ identifies with the
set  $\Hom\big(\pi_1(\Sigma),G\big)/G$ 
of equivalence classes of representations 
$\pi_1(\Sigma)\longrightarrow G$,
where two representations $\rho,\rho'$ are equivalent 
if and only if $\exists g\in G$ such that $\rho' = \Inn(g)\circ \rho$,
where $\Inn(g) : x \longmapsto gxg^{-1}$ is the inner automorphism
associated to $g\in G$. 
Since $\pi_1(\Sigma)$ is finitely generated, 
$\Hom\big(\pi_1(\Sigma),G\big)$ 
has the structure of a real-analytic subset in 
a Cartesian power $G^N$, and this structure is independent
of the choice of generators.
Give $\Hom\big(\pi_1(\Sigma),G\big)$ the classical topology
and note that it is  stratified into smooth submanifolds.
Give  $\Hom\big(\pi_1(\Sigma),G\big)/G$ the quotient topology.

The space $\Hom\big(\pi_1(\Sigma),G\big)/G$ may enjoy
several pathologies:
\begin{itemize}
\item The analytic variety 
$\Hom\big(\pi_1(\Sigma),G\big)$ may have singularities, and not
be a manifold;
\item $G$ may not act freely, even on the smooth points,
so the quotient map may be nontrivially branched,
and $\Hom\big(\pi_1(\Sigma),G\big)/G$ may have orbifold
singularities;
\item $G$ may not act properly, and the quotient
space
$\Hom\big(\pi_1(\Sigma),G\big)/G$ may not be Hausdorff.
\end{itemize}
All three pathologies may occur.

The automorphism group $\Aut\big(\pi_1(\Sigma)\big)$ 
acts on $\Hom\big(\pi_1(\Sigma),G\big)$ by right-composition.
The action of its subgroup $\Inn\big(\pi_1(\Sigma)\big)$ is 
absorbed in the $\Inn(G)$-action, and therefore the quotient
group
\begin{equation*}
\Out\big(\pi_1(\Sigma)\big) \;:=\;
\Aut\big(\pi_1(\Sigma)\big)/\Inn\big(\pi_1(\Sigma)\big)
\end{equation*}
acts on the quotient
\begin{equation*}
\Hom\big(\pi_1(\Sigma),G\big)/G.
\end{equation*}

Associating to a marked $(G,X)$-structure
the equivalence class of its holonomy representation
defines a continuous map
\begin{equation}\label{eq:hol}
\DGXS \xrightarrow{\hol} 
\Hom\big(\pi_1(\Sigma),G\big)/G
\end{equation}
which is evidently 
$\pi_0\big(\mathsf{Diff}(\Sigma)\big)$-equivariant, 
with respect to the homomorphism
\begin{equation*}
\pi_0\big(\mathsf{Diff}(\Sigma)\big)
\longrightarrow \Out\big(\pi_1(\Sigma)\big).
\end{equation*}

\begin{theorem*}[Thurston]
With respect to the above topologies, the holonomy map
$\hol$  in \eqref{eq:hol} is a local homeomorphism.
\end{theorem*}

For hyperbolic structures on closed surfaces, which 
are special cases of $(G,G)$-structures (or discrete embeddings
in Lie groups as above), this result is due to 
Weil~\cite{Weil_discretesubgroups, Weil2, Weil3}; 
see the very readable
paper by Bergeron-Gelander~\cite{BergeronGelander}.
This result is due to Hejhal~\cite{Hejhal} 
for $\cpo$-surfaces.
The general theorem  was first stated explicitly by Thurston~\cite{Thurston_1979notes}, 
and perhaps the first careful proof may be found 
in Lok~\cite{Lok} and  Canary-Epstein-Green~\cite{CEG}.
Bergeron and Gelander refer to this result as
the ``Ehresmann-Thurston theorem'' since many of the ideas
are implicit in Ehresmann's viewpoint~\cite{Ehresmann2}. 

The following  proof was worked out in \cite{geost}  with Hirsch, 
and was also known to Haefliger. By the covering homotopy theorem
and the local contractibility of $\Hom\big(\pi_1(\Sigma),G\big)$,
the isomorphism type of $E$ as a $G$-bundle is constant.
Thus one may assume that $E$ is a fixed $G$-bundle, 
although the flat structure (given by the transverse foliation $\Ff$)
varies, as the representation varies.
However it varies continuously in the $C^1$ topology.
Thus a given $\Ff$-transverse section $\Dev$ remains transverse
as $\Ff$ varies, and defines a geometric structure. 
This proves local surjectivity of $\hol$. 

Conversely, if $\Dev'$ is a transverse section sufficiently close to $\Dev$ 
in the $C^1$-topology, then it stays within a neighborhood of $\Dev(\Sigma)$. For a sufficiently small neighborhood $W$ of $\Dev(\Sigma)$, the foliation $\Ff|_W$ identifies with a product foliation of
$W\approx \Dev(\Sigma)\times X$ defined by the projection to $X$.
For each $m\in\Sigma$, the leaf of $\Ff|_W$ through $\Dev(m)$ 
meets $\Dev'(\Sigma)$ in a unique point $\Dev'(m')$ for $m'\in\Sigma$.
The correspondence $m\longmapsto m'$ is the required isotopy,
from which follows $\hol$ is locally injective.

\section{Thurston's geometrization of 3-manifolds}
In 1976, Thurston proposed that every closed 3-manifold admits a canonical decomposition into 
pieces, by cutting along surfaces of nonnegative Euler characteristic. Each of these pieces 
has one of eight geometries, modeled on eight $3$-dimensional Riemannian homogeneous spaces:
\begin{itemize}
\item {\bf {Elliptic geometry:\/}} Here $X = S^3$ and $G = \mathsf{O}(3)$ its group of isometries.
Manifolds with these geometries are the Riemannian $3$-manifolds of constant positive curvature, that is, {\em spherical space forms,\/}
and include lens spaces.
\item  $S^2\times\R$: The only closed $3$-manifolds with this geometry are $S^2\times S^1$ and a few quotients.
\item {\bf  {Euclidean geometry:\/}} Here $X = \R^3$ and $G$ its group of isometries.
These are the Riemannian manifolds of zero curvature,
and are quotients by torsionfree Euclidean {\em crystallographic
groups.\/} 
In 1912, Bieberbach proved every closed Euclidean manifold is a quotient of a flat torus by a finite group of isometries. 
Furthermore he proved there are only finitely many topological types of these manifolds, and that any homotopy-equivalence is homotopic
to an {\em affine isomorphism.\/}
\item {\bf  {Nilgeometry:\/}} Here again $X = \R^3$, regarded as the Heisenberg group with a left-invariant metric
and $G$ its group of isometries. Manifolds with these geometry are covered by nontrivial oriented $S^1$-bundles
over $2$-tori.
\item {\bf {Solvgeometry:\/}} Once again $X = \R^3$, regarded as a $3$-dimensional exponential solvable unimodular 
non-nilpotent Lie group and $G$ the group of isometries of a left-invariant metric. Hyperbolic torus bundles
(suspensions of Anosov diffeomorphisms of tori) have these structures.
\item $H^2\times\R$: Products of hyperbolic surfaces with $S^1$ have this geometry.
\item {\bf {Unit tangent bundle of $H^2$:\/}} An equivalent model is 
$\mathsf{PSL}(2,\R)$ with a left-invariant metric. Nontrivial oriented $S^1$-bundles of hyperbolic surfaces (such as the unit tangent bundle)
admit such structures.
\item {\bf {Hyperbolic geometry:\/}} Here $X = H^3$ and $G$ its group of isometries.
\end{itemize}
For a description of the eight homogeneous Riemannian geometries
and their relationship to $3$-manifolds, see the excellent surveys
by Scott~\cite{Scott} and Bonahon~\cite{Bonahon}.

\section{Complete affine 3-manifolds}

Manifolds modeled on {\em Euclidean geometry\/} are 
exactly the flat Riemannian manifolds. 
Compact Euclidean manifolds $M^n$
are precisely the quotients $\R^n/\Gamma$, where 
$\Gamma$ is a lattice of Euclidean isometries.
By the work of Bieberbach (1912), such a $\Gamma$ is
a finite extension of a lattice $\Lambda$ of translations. 
Thus $M$ is finitely covered by the torus $\R^n/\Lambda$.
Since all lattices $\Lambda\subset\R^n$ are {\em affinely\/}
the homotopy type of $M$ determines its affine equivalence class.
When $M$ is noncompact, but geodesically complete, then
$M$ is isometric to a flat orthogonal vector bundle over a compact Euclidean manifold.

These theorems give at least a qualitative classification of
manifolds with Euclidean structures. The generalization to
manifolds with {\em affine structures\/} is much more mysterious
and difficult. We begin by restricting to ones which are
{\em geodesically complete.\/} In that case the manifolds
are quotients $\R^n/\Gamma$ but $\Gamma$ is only assumed
to consist of affine transformations. 
However, unlike Euclidean manifolds considered above,
discreteness of $\Gamma\subset\mathsf{Aff}(\R^n)$ does not
generally imply the properness of the action, and the quotient
may not be Hausdorff. Characterizing which affine representations
define proper actions is a fundamental and challenging problem.

In the early 1960's, L.\ Auslander announced that every compact
complete affine manifold has virtually polycyclic fundamental
group, but his proof was flawed. 
In this case, the manifold is finitely covered by an
{\em affine solvmanifold\/} $\Gamma\backslash G$ where $G$ is
a (necessarily solvable) Lie group with a left-invariant complete
affine structure and $\Gamma\subset G$ is a lattice. 
Despite many partial results,
(\cite{FriedGoldman,AbelsMargulisSoifer,
AbelsMargulisSoifer2, Tomanov,GoldmanKamishima})
the {\em Auslander Conjecture\/} remains open.

Milnor~\cite{Milnor} asked whether the virtual polycyclicity of $\Gamma$
might hold even if the quotient $\R^n/\Gamma$ is {\em noncompact.\/}
Using the Tits Alternative~\cite{Tits}, he reduced this question
to whether a rank two free group could act {\em properly\/} by
affine transformations on $\R^n$. 
Margulis~\cite{Margulis
} showed, surprisingly, that such actions {\em do\/} exist when $n=3$.

For $n=3$, 
Fried and Goldman~\cite{FriedGoldman} showed that either $\Gamma$
is virtually polycyclic (in which case all the structures are easily
classified), or the linear holonomy homomorphism 
$\Gamma\xrightarrow{\mathsf{L}}\mathsf{GL}(3,\R)$ 
maps $\Gamma$ isomorphically onto a discrete  subgroup
of a conjugate of $\Oto\subset\mathsf{GL}(3,\R)$.
Since $\mathsf{L}^{-1}\Oto$ preserves a flat Lorentz metric
on $\R^3$, the geometric structure on $M$ refines to a 
flat Lorentz structure, modeled on $\E^3_1$, which is $\R^3$
with the corresponding flat Lorentz metric.
In particular $M^3 = \E^3_1/\Gamma$
is a {\em complete flat Lorentz 3-manifold\/}
and $\Sigma:= H^2/\mathsf{L}(\Gamma)$ is a complete hyperbolic surface. This establishes the Auslander Conjecture in dimension 3: 
the cohomological dimension of $\Gamma\cong\pi_1(M^3)$ 
equals $3$ since $M$ is aspherical, but the cohomological dimension
$\Gamma\cong\pi_1(\Sigma)$  is at most $2$. 
In 1990, Mess~\cite{Mess} proved that the surface $\Sigma$ is 
{\em noncompact,\/} and therefore $\Gamma$ must be a free group.
(Compare also Goldman-Margulis~\cite{GoldmanMargulis} and Labourie~\cite{Labourie_fuchsian} for other proofs.)

Drumm~\cite{Drumm0,Drumm2} (see also \cite{CharetteGoldman} )
gave a geometric construction of these quotient 
manifolds using polyhedra in Minkowski space $\R^3_1$ 
now called {\em crooked planes.\/} Using crooked planes, 
he showed that every noncompact complete hyperbolic surface
$\Sigma$ arises from a complete flat Lorentz $3$-manifold; 
that is, he showed that every non-cocompact Fuchsian group
$\mathsf{L}(\Gamma)\subset\Oto$ admits a {\em proper\/}
affine deformation $\Gamma$.

The conjectural picture of these manifolds is as follows.

The space of equivalence classes of affine deformations of $\Gamma$
is the vector space $H^1(\Gamma,\R^3_1)$, 
and the proper affine deformations define an open convex
cone in this vector space. 
Goldman-Labourie-Margulis~\cite{GLM}
have proved this when $\Gamma$
is finitely generated  and contains no parabolic elements.
Furthermore a finite-index subgroup of $\Gamma$ should have a fundamental domain which is bounded by crooked planes,
and $M^3$ should be homeomorphic to a solid handlebody.
Charette-Drumm-Goldman~\cite{CharetteDrummGoldman} 
have proved  this when $\Sigma$ is homeomorphic to a $3$-holed sphere.

Translational conjugacy classes of affine deformations of a Fuchsian
group $\Gamma_0\subset\Oto$ comprise the cohomology group
$H^1(\Gamma_0;\R^3_1)$. As the $\Oto$-module $\R^3_1$ identifies
with the Lie algebra of $\Oto$ with the adjoint representation, 
this cohomology group identifies with the {\em space of infinitesimal
deformations of the hyperbolic surface\/} $\Sigma = H^2/\Gamma_0$.
(Compare Goldman-Margulis~\cite{GoldmanMargulis} and \cite{Goldman_MargulisInvariant}.)

When $\Sigma$ has no cusps, \cite{GLM}
provides a criterion for properness of an affine deformation corresponding to a deformation $\sigma$ of the hyperbolic surface $\Sigma$.
The affine deformation $\Gamma_\sigma$ acts properly
on $\E^3_1$ if and only
if every probability measure on $U\Sigma$ invariant under the 
geodesic flow {\em infinitesimally lengthens\/} (respectively
{\em infinitesimally shortens\/} under $\sigma$. (We conjecture
a similar statement in general.) 
Using ideas based on Thurston~\cite{stretch}, one can
reduce this to probability measures arising from measured
geodesic laminations. When $\Sigma$ is a three-holed sphere,
\cite{CharetteDrummGoldman} 
implies the proper affine deformations are precisely the ones for which the three components of $\partial\Sigma$ either all infinitesimally lengthen or all infinitesimally shorten.

Other examples of {\em conformally flat Lorentzian manifolds\/}
have recently been studied by Frances~\cite{Frances}, Zeghib~\cite{Zeghib}, and Bonsante-Schlenker~\cite{BonsanteSchlenker}, 
also closely relating to hyperbolic geometry.

\section{Affine structures on closed manifolds}

The question of which closed manifolds admit affine structures
seems quite difficult. Even for complete structures, the pattern
is mysterious.  Milnor~\cite{Milnor} asked whether every 
virtually polycyclic group arises as the fundamental group of 
a compact complete affine manifold.
Benoist~\cite{Benoist_nilvarietes,
Benoist_nonaffine} found $11$-dimensional nilpotent counterexamples.
However by replacing $\R^n$ by a simply connected nilpotent Lie group,
one obtains more general structures.  
Dekimpe~\cite{Dekimpe} showed that every virtually polycyclic
group arises as the fundamental group of such a {\em NIL-affine\/}
manifold.

For incomplete structures, 
the picture is even more unclear.  
The {\em Markus conjecture, \/} first stated by L.\ Markus as a homework exercise in unpublished lecture notes at the University of Minnesota in 1960 asserts that, for closed affine manifolds, geodesic completeness
is equivalent to parallel volume (linear holonomy in $\mathsf{SL}(n,\R)$.
That this conjecture remains open testifies to our current ignorance.

An important partial result is Carri\`ere's result~\cite{Carriere}
that a closed flat Lorentzian manifold is geodesically complete.
This has been generalized in a different direction by Klingler~\cite{Klingler} to all closed Lorentzian manifolds with {\em constant\/}
curvature.

Using parallel volume forms,
Smillie~\cite{Sm3} showed that the holonomy of a compact affine manifold cannot factor through a free product of finite groups; 
his methods were extended by Goldman-Hirsch~\cite{GH3,GH4}
to prove nonexistence results for affine structures on closed manifolds
with certain conditions on the holonomy. Using these results,
Carri\`ere, Dal'bo and Meigniez~\cite{CDM} showed that certain Seifert 3-manifolds with hyperbolic base admit no affine structures.

Perhaps the most famous conjecture about affine structures on closed manifolds is Chern's conjecture that a closed affine manifold must have
Euler characteristic zero. For flat pseudo-Riemannian manifolds or 
{\em complex\/} affine manifolds, this follows from Chern-Gauss-Bonnet.
Using an elegant argument, Kostant and Sullivan~\cite{KostantSullivan}
proved this conjecture for complete affine manifolds. (This would follow
immediately from the Auslander Conjecture.)

In a different direction, Smillie~\cite{Sm1} found simple examples of
closed manifolds with flat tangent bundles (these would have affine connections with zero curvature, but possibly nonzero torsion). 
Recent results in this direction have been obtained by Bucher-Gelander~\cite{BucherGelander}.

\section{Hyperbolic geometry on 2-manifolds}
The prototype of geometric structures, and historically
one of the basic examples, are hyperbolic structures
on surfaces $\Sigma$ with $\chi(\Sigma)<0$. 
Here $X$ is the hyperbolic plane and 
$G\cong\mathsf{PGL}(2,\R)$. 
Fricke and Klein~\cite{FrickeKlein} studied the deformation space of hyperbolic structures on $\Sigma$ as well as on $2$-dimensional orbifolds. The deformation space $\Fricke$ 
of marked hyperbolic structures on $\Sigma$
(sometimes called {\em Fricke Space\/} (\cite{BersGardiner})
can also be described as the space of equivalence classes
of discrete embeddings $\pi_1(\Sigma)\longrightarrow G$.
The  Poincar\'e-Klein-Koebe Uniformization
Theorem relates hyperbolic structures and complex structures,
so the Fricke space identifies with the {\em Teichm\"uller space\/}
of $\Sigma$, which parametrizes Riemann surfaces homeomorphic
to $\Sigma$. For this reason, although Teichm\"uller himself 
never studied hyperbolic geometry, the deformation theory of hyperbolic
structures on surfaces is often referred to as {\em Teichm\"uller theory.}

Representations of surface groups in $G=\PSLtR$ closely
relate to geometric structures. 
A representation $\pi_1(\Sigma)\xrightarrow{\rho}G$
determines an oriented flat $H^2$-bundle over $\Sigma$.
Oriented flat $H^2$-bundles are classified by their {\em Euler class\/}, which lives in $H^2(\Sigma;\Z)\cong \Z $ when $\Sigma$ is closed and oriented.
The Euler number of a flat oriented $H^2$-bundle satisfies
\begin{equation}\label{eq:MilnorWood}
\vert \Euler(\rho)\vert \le  -\chi(\Sigma)
\end{equation}
as proved by Wood\cite{Wood}, following earlier work of 
Milnor\cite{Milnor}.

\begin{theorem}\label{thm:Goldman_thesis}
Equality holds in \eqref{eq:MilnorWood} if and only if
$\rho$ is a discrete embedding.
\end{theorem}
This theorem was first proved in \cite{Goldman_thesis},
using Ehresmann's viewpoint. Namely, the condition that
$\Euler(\rho) = \pm \chi(\Sigma)$ means that the associated
flat $H^2$-bundle $E_\rho$ with holonomy  homomorphism $\rho$ 
is isomorphic (up to changing orientation)
to the tangent bundle of $\Sigma$ (as a topological disc bundle, or equivalently a microbundle over $\Sigma$).
If $\rho$ is the holonomy of a hyperbolic surface $M\approx \Sigma$,
then $E(M) = E_\rho\approx T\Sigma$.  
Theorem~\ref{thm:Goldman_thesis} is a converse:
if the flat bundle ``is isomorphic to the tangent bundle
(as a $(G,X)$-bundle)",  then 
the flat $(G,X)$-bundle arises from a $(G,X)$-structure on $\Sigma$.

In the case the representation $\rho$ has discrete torsionfree cocompact image, Theorem~\ref{thm:Goldman_thesis} reduces to a classical result of Kneser~\cite{Kneser}. In 1930 Kneser proved that if $\Sigma\xrightarrow{f}\Sigma'$ is a continuous map of degree $d$, then 
\begin{equation*}
d \vert \chi(\Sigma') \vert \;\le \; 
\vert \chi(\Sigma) \vert 
\end{equation*}
with equality $\Longleftrightarrow$ $f$ is homotopic to a covering space.
(In this case  $\Sigma'$ is the hyperbolic surface obtained as the 
quotient by the image of $\rho$, and $\Euler(\rho) = d  \chi(\Sigma') $.
Kneser's theorem is thus a {\em discrete\/} version of 
Theorem~\ref{thm:Goldman_thesis}.

By now Theorem~\ref{thm:Goldman_thesis}
has many proofs and extensions.
One proof, using harmonic maps,
begins by choosing a Riemann surface $M\approx \Sigma$.
Then, by Corlette~\cite{Corlette} and Donaldson~\cite{Donaldson},
either the image of $\rho$ is solvable (in which case $\Euler(\rho) = 0$)
or the image is reductive, and there exists a $\rho$-equivariant 
harmonic map $\tM \xrightarrow{h} X$. By an adaptation of 
Eels-Wood~\cite{EelsWood}, $\Euler(\rho)$ can be computed
as the sum of local indices of the critical points of  $h$.
In particular, the assumption of {\em maximality:\/} 
$\Euler(\rho) = \pm \chi(M)$ implies that $h$ must be
holomorphic (or anti-holomorphic), and using the arguments
of Schoen-Yau~\cite{SchoenYau_harmonic}, $h$ must be a diffeomorphism.
In particular $\rho$ must be a discrete embedding.

Shortly after \cite{Goldman_thesis}, another proof was given by Matsumoto~\cite{Matsumoto} (compare also Mess~\cite{Mess}),
related to ideas of bounded cohomology. 
This led to the work of Ghys~\cite{Ghys}, who proved that
the Euler class of an orientation-preserving action of $\pi_1(\Sigma)$ on $S^1$ is a {\em bounded\/} class, and its class in bounded 
cohomology determines the action up to topological semi-conjugacy.
In particular maximality in the Milnor-Wood inequality \eqref{eq:MilnorWood} implies the topological action is 
conjugate to the projective action arising from (any)  
discrete embedding in $\PSLtR$.

The Euler number classifies components of
$\Hom\big(\pi_1(\Sigma),\PSLtR)$.
That is, if $\Sigma$ is closed,
oriented, of genus $g>1$, the $4g-3$ connected components
are the inverse images $\Euler^{-1}(j)$ where 
\begin{equation*}
j\; =\; 2-2g,3-2g,\dots, 2g-2
\end{equation*}
(Goldman~\cite{Goldman_components}).
Independently, Hitchin~\cite{Hitchin_selfduality} gave a much 
different proof, using Higgs bundles. 
Moreover he identified the Euler class $2-2g + k$ component
with a vector bundle over the $k$-th symmetric power of $\Sigma$
(compare the expository article ~\cite{Goldman_Hitchin})

When $G$ is a semisimple {\em compact\/} or {\em complex\/}
Lie group, components of the representation space bijectively
correspond to  $\pi_1(G)$. In particular in these basic cases, 
the number of components  is {\em independent\/} of the genus. 
(See Li~\cite{Li} and Rapinchuk--Benyash-Krivetz--Chernousov\cite{RapinchukBKC}.) Recently Florentino and Lawton~\cite{FlorentinoLawton} have determined the homotopy type
of $\Hom(\Gamma,G)//G$ when $\Gamma$ is free and $G$ is a complex reductive group.

This simple picture becomes much more intricate
and fascinating for higher dimensional noncompact real Lie groups; 
the most effective technique so far has been the interpretation
in terms of Higgs bundles and the use of infinite-dimensional 
Morse theory; see Bradlow-Garcia-Prada-Gothen~\cite{BradlowGarciaPradaGothen} for a survey of some recent
results on the components when $G$ is a simple {\em real\/} Lie group.

Theorem~\ref{thm:Goldman_thesis} leads to rigidity theorems
for surface group representations as well. When $G$ is the automorphism group of a Hermitian symmetric space $X$, integrating a $G$-invariant K\"ahler form on $X$ over a smooth section of
a flat $(G,X)$-bundle induces a characteristic class 
$\tau(\rho)$ 
first defined by Turaev~\cite{Turaev} and Toledo~\cite{Toledo_89}. 
This characteristic class satisfies an inequality similar to
\eqref{eq:MilnorWood}. 
The {\em maximal representations,\/} (when equality is attained)
have very special properties.
When $X$ is complex hyperbolic space,  a representation $\pi_1(\Sigma)\xrightarrow{\rho}\PUno$ is maximal if and only if 
it stabilizes a totally geodesic holomorphic curve, and its
restriction  is  Fuchsian (Toledo~\cite{Toledo_89}).

In higher rank the situation is much more
interesting and complicated. 
Burger-Iozzi-Wienhard~\cite{BurgerIozziWienhard_toledo}
showed that maximal representations are discrete embeddings,
with reductive Zariski closures. With Labourie, 
they proved~\cite{Burger_Iozzi_Labourie_Wienhard} in the case of
$\Spnr$, that these representations quasi-isometrically embed
$\pi_1(\Sigma)$ in $G$. Many of these properties follow from
the fact that maximal representations are Anosov representations
in the sense of Labourie~\cite{Labourie_anosov}.  
Using  Higgs bundle theory, Bradlow-Garcia-Prada-Gothen~\cite{BradlowGarciaPradaGothen} have counted components
of maximal representations.  Guichard-Wienhard~\cite{GuichardWienhard2} have found components of maximal
representations in $\Spnr$, all of whose elements have Zariski
dense image (in contrast to $\PUno$ discussed above).
For a good survey of these results, see
Burger-Iozzi-Wienhard~\cite{BurgerIozziWienhard_handbook}.

\section{Complex projective 1-manifolds, flat conformal
structures and spherical CR structures}

When $X$ is enlarged to $\cpo$ and $G$ to $\PSLtC$,
the resulting deformation theory of $\cpo$-structures is quite
rich. A manifold modeled on this geometry is naturally a Riemann
surface, and thus the deformation space fibers over the Teichm\"uller
space of marked Riemann surfaces:
\begin{equation}\label{eq:poincare}
\DGXS \longrightarrow \Teich.
\end{equation}
The classical theory of the Schwarzian derivative 
identifies this fibration with a holomorphic affine bundle,
where the fiber over a point in $\Teich$ corresponding to
a marked Riemann surface $\Sigma \xrightarrow{\approx} M$
is an affine space with underlying vector space $H^0(M;\kappa_M^2)$
consisting of holomorphic quadratic differentials on $M$.

In the late 1970's, Thurston (unpublished) showed that 
$\DGXS$ admits an alternate description as $\Fricke\times\mathcal{ML}(\Sigma)$ where $\mathcal{ML}(\Sigma)$ is the space of equivalence classes of measured geodesic laminations on $\Sigma$. 
(Compare Kamishima-Tan~\cite{KamishimaTan}.)
\cite{fuchsianholonomy} gives the topological classification of $\cpo$-structures whose holonomy representation is a quasi-Fuchsian
embedding.
Gallo-Kapovich-Marden~\cite{GalloKapovichMarden} showed
that the image of the holonomy map $\hol$ consists of representations
into $\PSLtC$ which lift to an irreducible and unbounded representation into $\SLtC$.

For an excellent survey of this subject, see Dumas~\cite{Dumas}.

These structures generalize to higher dimensions in several
ways. For example $\PSLtC$ is the group of orientation-preserving conformal automorphisms of $\cpo\approx S^2$.
A {\em flat conformal structure\/} is a geometric structure locally
modeled on $S^n$ with its group of conformal automorphisms.
This structure is equivalent to a conformal class of Riemannian
metrics, which are {\em locally\/}  conformally equivalent to
Euclidean metrics.  (Compare Matsumoto~\cite{Matsumoto_conf}.)
In the 1970's it seemed tempting to try to prove the 
Poincar\'e conjecture by showing that every closed $3$-manifold
admits such a structure. This was supported by the fact that
these structures are closed under connected sums (Kulkarni~\cite{Kulkarni}). 
This approach was further promoted by the
fact that such structures arise as critical points of the Chern-Simons
functional ~\cite{ChernSimons}, 
and one could try to reach critical points by following the gradient
flow of the Chern-Simons functional.  However,
closed $3$-manifolds with nilgeometry or solvgeometry admit no
flat conformal structures whatsoever~\cite{conf}). 

As $H^{n-1}\times \R$ embeds in $S^n$ as the complement
of a codimension-two subsphere, 
the conformal geometry of $S^n$ contains 
$H^{n-1}\times \R$-geometry. Thus products of closed
surfaces with $S^1$ do admit flat conformal structures, and
Kapovich~\cite{Kapovich} and 
Gromov-Lawson-Thurston~\cite{GromovLawsonThurston} 
showed that even some nontrivial $S^1$-bundles over closed
surfaces admit flat conformal structures,
although $T_1(H^2)$-geometry admits no conformal model in $S^3$.

Kulkarni-Pinkall~\cite{KulkarniPinkall} have extended
Thurston's correspondence 
\begin{equation*}
\DGXS \longleftrightarrow \Fricke\times\mathcal{ML}(\Sigma)
\end{equation*}
to associate to a flat conformal structure on a manifold 
(satisfying a generic condition of ``hyperbolic type'') 
a hyperbolic metric with some extrinsic (bending) data.
b

A similar class of structures are the {\em spherical CR-structures,\/}
modeled on $S^{2n-1}$ as the boundary of {\em complex hyperbolic 
$n$-space,\/}  in the same way that $S^{n-1}$ with its conformal structure bounds  real hyperbolic $n$-space. Some of the first
examples were given by Burns-Shnider~\cite{BurnsShnider}.
$3$-manifolds with nilgeometry naturally admits such structures, 
but by~\cite{conf}, closed $3$-manifolds with Euclidean and solvgeometry do not admit such structures. 
Twisted $S^1$-bundles admit many such structures
(see for example \cite{GoldmanKapovichLeeb}), but 
recently Ananin, Grossi and Gusevskii~\cite
{AGG,AG} %
have constructed surprising examples of spherical CR-structures
on products of closed hyperbolic surfaces with $S^1$. 
Other interesting examples of spherical CR-structures on 
$3$-manifolds have been constructed by Schwartz~\cite{Schwartz_ICM,Schwartz_Dehn,Schwartz_GoldmanParker}, 
Falbel~\cite{Falbel}, Gusevskii,  Parker~\cite{Parker_Lattices},
Parker-Platis~\cite{ParkerPlatis}.

When $X = \rpn$ and $G= \PGLnoR$,
then a $(G,X)$-structure is a {\em flat projective connection.}

In dimension $3$, the only closed manifold known {\em not\/}
to admit an $\rpthree$-structure is the connected sum 
$\rpthree \# \rpthree$ (Cooper-Goldman~\cite{CooperGoldman}).
Many diverse examples of $\rpthree$-structures on twisted $S^1$-bundles over closed hyperbolic surfaces arise from maximal
representations of surface groups into $\Spfour$ 
by Guichard-Wienhard~\cite{GuichardWienhard2}.
All eight of the Thurston geometries have models in $\rpthree$~\cite{Molnar}.

The $2$-dimensional theory is relatively mature.
The most important examples are the {\em convex\/}
structures, namely those which arise as quotients $\Omega/\Gamma$
where $\Omega$ is a convex domain in $\rpt$
and $\Gamma$ is a group of collineations preserving $\Omega$. Kuiper~\cite{Kuiper_affine} showed that all
convex  structures on $2$-tori are affine structures, and classified them.
They are all quotients of the plane, a half-plane or a quadrant.
In higher genus, he showed~\cite{Kuiper_rpt} 
that either $\partial\Omega$ is a conic
(in which case the projective structure is a hyperbolic structure)
or it fails to be $C^2$.  Benzecri~\cite{Benzecri2} showed that in the
latter case, it is $C^1$ and is strictly convex. Using the analog
of Fenchel-Nielsen coordinates, Goldman~\cite{Goldman_convex}
showed that the deformation space $\Cc$ is a cell of dimension
$-8\chi(\Sigma)$.  (Kim~\cite{HongChanKim} showed these coordinates
are global Darboux coordinates for the symplectic structure, 
extending a result of Wolpert~\cite{Wolpert_curves} for $\Fricke$.)
In his doctoral thesis, Choi showed that every
structure on a closed surface {\em canonically\/} decomposes
into convex structures with geodesic boundary, glued together
along boundary components. Combining these two results,
one identifies the deformation space {\em precisely\/} as
a countable disjoint union of open 
$-8\chi(\Sigma)$-cells~\cite{Goldman_Choi2}.

Using analytic techniques, Labourie~\cite{Labourie_margulis}
and Loftin~\cite{Loftin}, independently, described $\Cc$ as a cell in a quite different way. Associated to a convex $\rpt$-structure $M$ is
a natural Riemannian metric arising from representing $M$
as a convex surface in $\R^3$, which is a {\em hyperbolic affine sphere.\/} The underlying conformal structure defines a point 
in $\Teich$ associated to the convex $\rpt$-manifold $M$. 
Its extrinsic geometry is described by a {\em holomorphic cubic differential\/} on the corresponding Riemann surface. 
In this way $\Cc$ identifies with the bundle over $\Teich$
whose fiber over a marked Riemann surface is the vector space of
holomorphic cubic differentials on that Riemann surface.
Loftin~\cite{Loftin2} relates  the geometry of these structures to the asymptotics of this deformation space.

These results generalize in several directions. 
In a series of beautiful papers, 
Benoist~\cite{Benoist_convex,
Benoist_1,Benoist_2,Benoist_3,Benoist_4,
Benoist_survey} 
studied convex projective structures $\Omega/\Gamma$ on compact manifolds.
The natural {\em Hilbert metric\/} on $\Omega$ determines a 
(Finsler) metric on $M$, and if $\Omega$ is strictly convex,
then this natural metric has negative curvature and $\Gamma$
is a hyperbolic group. The corresponding geodesic flow is an
Anosov flow, which if $M$ admits a hyperbolic structure, is topologically conjugate to the geodesic flow of the hyperbolic metric.
Furthermore, as in \cite{Goldman_Choi1}, the corresponding
representations $\Gamma\longrightarrow\mathsf{PGL}(n+1,\R)$
form a connected component of the space of representations.
For compact quotients $\Omega/Gamma$, Benoist showed that the hyperbolicity of the group $\Gamma$ is equivalent to the strict
convexity of $\partial\Omega$. He constructed $3$-dimensional 
examples of convex structures on $3$-manifolds with incompressible
tori and hyperbolic components, where $\partial\Omega$ is the 
closure of a disjoint countable union of triangles. 
In a different direction, Kapovich~\cite{Kapovich_GT} constructed
convex projective structures with $\partial\Omega$ strictly convex
but $\Omega/\Gamma$ has no locally symmetric structure.


When $G$ is a split real form of a complex semisimple Lie group,
Hitchin~\cite{Hitchin_Lie} showed that $\HompGG$ contains
components homeomorphic to open cells. 
Specifically, these are the components containing Fuchsian
representations into $\SLtR$ composed with the Kostant
principal  representation $\SLtR\longrightarrow G$. 
When $G=\SLthreeR$, then $\hol$ maps $\Cc$ diffeomorphically
to Hitchin's component (Choi-Goldman~\cite{Goldman_Choi1}).
Guichard and Wienhard~\cite{GuichardWienhard1}
have found interpretations of Hitchin components in $\SLfourR$
in terms of geometric structures. Recently~\cite{GuichardWienhard3} they have also shown
that a very wide class of {\em Anosov representations\/} as
defined by Labourie~\cite{Labourie_anosov},  correspond to
geometric structures on {\em closed manifolds.\/}
(A much different class of Anosov representations of surface
groups has recently been studied by Barbot~\cite{Barbot1,Barbot2}.

The properness of the action of $\Mod$ on $\Fricke$ is generally
attributed to Fricke. Many cases are known of components
of deformation spaces when $\Mod$ acts 
properly~\cite{Goldman_Wentworth,Wienhard_mapping,
Burger_Iozzi_Labourie_Wienhard}. In many of these cases, these components consist of holonomy
representations of uniformizable Ehresmann structures.

\section{Surface groups: symplectic geometry and mapping class group}

Clearly the classification of geometric structures in low dimensions
closely interacts with the space of surface group representations.
Many examples have already been discussed here.
By the Ehresmann-Weil-Thurston  holonomy theorem,
the local geometry of 
$\HompGG$ 
is the same local geometry of $\DGXS$. When $\Sigma$
is a compact surface, this space itself admits rich geometric structures.

Associated to an orientation on $\Sigma$ and an $\Ad(G)$-invariant
nondegenerate symmetric bilinear form $\Bb$ on the Lie algebra of $G$
is a natural {\em symplectic structure\/} on the deformation space.
(When $\partial\Sigma\neq\emptyset$, one obtains a {\em Poisson 
structure\/} whose symplectic leaves correspond to fixing the
conjugacy classes of the holonomy along boundary components.)
This extends the cup-product symplectic structure on
$H^1(\Sigma,\R)$ (when $G=\R$), the K\"ahler {\em form\/}
on the Jacobian of a Riemann surface $M\approx \Sigma$,
(when $G=\mathsf{U}(1)$), and the Weil-Petersson K\"ahler
form on $\Teich$ (when $G=\PSLtR$). 
Compare \cite{Goldman_symp}.

The symplectic geometry extends over the singularities of
the deformation space as well. 
In joint work with Millson~\cite{Goldman_Millson2,Millson_ICM}, 
inspired by a letter of Deligne~\cite{Deligne},
it is shown that the germ at a reductive representation $\rho$, the 
analytic variety $\HompG$ is locally equivalent to a cone defined by
a system of homogeneous quadratic equations.
Explicitly, this quadratic cone is defined by the cup-product
$$
\Zo\times\Zo  \xrightarrow{[,]_*\cup} 
\mathsf{H}^2\Sigma,\mathfrak{g}_{\Ad \rho})
$$
using Weil's identification of the Zariski tangent space of
$\HompG$ at $\rho$ with $\Zo$.
This quadratic singularity theorem extends to higher-dimensional
K\"ahler manifolds~\cite{Simpson} and relates to the 
stratified symplectic spaces considered by  Sjamaar-Lerman~\cite{SjamaarLerman}.

The symplectic/Poisson geometry of the deformation spaces 
$\HompGG$ and $\DGXS$ associate vector fields to functions 
in the following way (see \cite{Goldman:invariantfunctions}. 
A natural class of functions $f_\alpha$ 
on $\HompGG$ arise from $\Inn(G)$-invariant functions $G\xrightarrow{f}\R$ and elements $\alpha\in\pi(\Sigma)$ by composition:
$$
[\rho] \;\xrightarrow{f_\alpha}   f\big(\rho(\alpha)\big).
$$
For example, when $\ell$ is the geodesic length function on 
$\PSLtR$, this construction yields the geodesic length functions 
$\ell_\alpha$ on $\Teich$.

When $\alpha$ arises from a {\em simple closed curve\/} on $\Sigma$ 
then the Hamiltonian flow associated to the vector field $\Ham(f_\alpha)$
admits a simple description as a {\em generalized twist flow.\/}
Such a flow is ``supported on $\alpha$'' in the sense that pulled back
to the complement $\Sigma\setminus\alpha$ the flow is a trivial
deformation. This extends the results
of Wolpert~\cite{Wolpert_FenchelNielsen,Wolpert_symplectic} for 
the Weil-Petersson symplectic form on  $\Teich$, 
Fenchel-Nielsen twist flow (or {\em earthquake\/}) along $\alpha$ 
is $\Ham(\ell_\alpha)$. 
For the case of $G=\mathsf{SU}(2)$, Jeffrey and Weitsman~\cite{JeffreyWeitsman} used these flows to define an ``almost
toric'' structure on $\HompGG$ from which they deduced the
Verlinde formulas.

The Poisson brackets of the functions $f_\alpha$ may be computed in terms of oriented intersections on $\Sigma$. For $G=\mathsf{GL}(n)$, and $f=\tr$, one obtains a topologically defined Lie algebra based on homotopy classes of curves on $\Sigma$ with a representation in the Poisson algebra of functions on $\HompGG$. Turaev\cite{Turaev2} showed this Lie algebra extends to a Lie bialgebra and found several quantizations. Recently Moira Chas~\cite{Chas} has discovered algebraic properties of this Lie algebra; in particular she proved that 
the $\ell^1$ norm of a bracket $[\alpha,\beta]$ of two unoriented simple closed curves equals the geometric intersection number
 $i(\alpha,\beta)$.

These algebraic structures  extend in 
higher dimensions to the  {\em string topology\/}  of  Chas-Sullivan~\cite{ChasSullivan}.

The symplectic geometry is $\Mod$-invariant 
and in particular defines an invariant measure on
the deformation space. Unlike the many cases in which $\Mod$
acts properly discussed above, when $G$ is compact,
this measure-preserving action is ergodic on each connected
component (Goldman~\cite{Goldman_erg}, Pickrell-Xia~\cite{PickrellXia}, Goldman-Xia~\cite{GoldmanXia_Zimmer}).
When $G$ is noncompact, invariant open subsets of the deformation space exist where the action is proper (such as the subset of Anosov
representations), but in general $\Mod$ can act properly 
on open subsets containing non-discrete representations, even
for $\PSLtR$ (\cite{Goldman_modular,Goldman_McShane_Stantchev_Tan, TYZ}).

Similar questions for the action of the outer automorphism group
$\Out(\mathbb{F}_n)$ of a free group $\mathbb{F}_n$ 
on $\Hom(\mathbb{F}_n,G)/G$ have recently been 
studied~\cite{Goldman_outer}. In particular Gelander has proved
that the action of $\Out(\mathbb{F}_n)$ is ergodic whenever
$G$ is a compact connected Lie group. For $G=\SLtC$, 
Minsky~\cite{Minsky} has recently found open subsets of 
$\Hom(\mathbb{F}_n,G)/G$ strictly containing the subset of 
Schottky embeddings for which the action is proper.

\section*{Acknowledgement}
I would like to thank Virginie Charette, Son Lam Ho,
Aaron Magid,
Karin Melnick, 
and Anna Wienhard for helpful suggestions in the preparation of this
manuscript.



\end{document}